\documentclass[a4paper,10pt]{article}

\newtheorem{thm}{Theorem}

\newtheorem{rem}[thm]{Remark}
\newtheorem{lem}[thm]{Lemma}
\newtheorem{prop}[thm]{Proposition}
\def \a{{\alpha}}
\def \b{{\beta}}

\def \t{{\vartheta}}

\def \m{{\mu}}

\def \E{{\bf E}\, }

\def \P{{\bf P}}

\def \qq{{\qquad}}
\def \R{{\bf R}}

\def \noi{{\noindent}}

\def\qed{\hbox{\vrule height 6pt depth 0pt width
6pt}}
\def\cqfd{\hfill\penalty 500\kern 10pt\qed\medbreak}

\font\phh=cmcsc10
at  8 pt
\scrollmode
\title{A Remark on Zeros   of Brownian Motion}
 \author{
   Michel Weber}

\begin{document}

\maketitle
%%%%%%%%%%%%%%%%%%%%%%%%%%%%%%%%%%%%%%%%%%%%%%%%%%%%%%%%%%%%%%%%%%%%%%%%%%%%%%%%%%%%%%%%%%%%%%

\begin{abstract}
 Let
$ \{W(t), t\ge 0\}$ be a standard Brownian motion. If  $I$ is a bounded interval  on which   $W $ has no zero, an almost sure lower bound to  $\inf\{
|W(t)|, t\in I\}$ can be provided, when $I$  is taken from a given countable family of intervals covering the positive half-line.   
\end{abstract}

%%%%%%%%%%%%%%%%%%%%%%%%%%%%%%%%%%%%%%%%%%%%%%%%%%%%%%%%%%%%%%%%%%%%%%%%%%%%%%%%%%%%%%%%%%%%%%%

\section{Main Result}
  Let
$ \{W(t), t\ge 0\}$ be a standard Brownian motion. Let $I$ be some bounded interval of $\R^+$. Suppose $W(t)\not= 0$, for all $t\in I$. What can be
  said about the size of $\inf\{ |W(t)|, t\in I\}$? This one only depends on the location of $I$   and of the size of $I$. The object of
this note is to prove more  precisely the following result. 
\begin{thm}  Let $\t_k\ge 0$ be such that $T_N=\sum_{k\le N}\t_k\uparrow
\infty$ and denote $I_N=[T_{N }, T_{N+1}]$. Let $ \eta_k\ge 0$ be such that 
$$\sum_{N\ge 1}\eta_N \min\Big(   { 1  
\over  \sqrt{ T_{N+1}-T_{N } }}      \, ,{1\over \sqrt{   T_{N }} }\Big) <\infty$$
Then  
$$\P\Big\{ \inf_{t\in I_N}|W(t)|\ge  \eta_N\ {or}\   W(t) =0 \ {for \ some}\ t\in I_N, \quad   N\ { 
ultimately}\Big\}=1.
$$
\end{thm}The proof  relies upon several intermediate results on infima of $|W|$, which are also of independent interest. 
\section{Local infima of Brownian motion}
   In this section, we collect some  properties of the infimum of $W$ over bounded
intervals. Precise estimates of the probability 
$$\P\big\{\inf_{a\le t\le b}|W(t)-M|\ge  c \big\}.$$
 will be necessary. Notice preliminary, since $-W$ and $W$ have same law that
\begin{eqnarray*}\P\big\{\inf_{a\le t\le b}|W(t)+M|\ge  c \big\}&=&\P\big\{\inf_{a\le t\le b}|-W(t)+M|\ge  c \big\}\cr
&=&\P\big\{\inf_{a\le
t\le b}|W(t)- M |\ge  c \big\},
\end{eqnarray*}
so that it is enough to consider the case $M\ge 0$. Put 
$$\Psi(x)= \P\{ W(1)>x\}= \int_x^\infty e^{-u^2/2} {du \over \sqrt {2\pi}}, \qq \quad x\in \R.$$
The    lemma  below is certainly well-known, although we could
not find a reference. We included a proof for the sake of completeness.
\begin{lem} \label{l2}Let $0<a<b<\infty$. Then for any $c> 0$ and any real $M$  
$$ \P\big\{\inf_{a\le t\le b}|W(t)-M|\ge  c \big\} =\int_{|v|>c } \Big[1-2\Psi({|v|-c\over \sqrt{b-a}})\Big]{ e^{-{(M+v)^2\over
2a}}\over
\sqrt{2\pi a}}dv   . $$
\end{lem}
 \noi   {\it Proof.}   By symmetry of the law of  $W$    it suffices to consider the case $M\ge 0$.
  By the intermediate values Theorem,  
\begin{equation}\P\big\{\inf_{a\le t\le b}|W(t)-M|\ge c \big\}= \P\big\{\inf_{a\le t\le b} W(t) \ge  M+c\big\}+\P\big\{ \sup_{a\le t\le b}
W(t) \le M-c\big\}. 
\end{equation}

   Let $x\ge 0$. Then $ \P\big\{\inf_{a\le t\le b}|W(t)-M|\ge  c\,\big|\,W(a) =M\pm
x\big\}=0$, if
$0\le x\le c$; and if $x>c$, 
$$ \P\big\{\inf_{a\le t\le b}|W(t)-M|\ge  c\,\big|\,W(a) =M+x\big\}=\P\big\{\sup_{a\le t\le b}
(W(a)-W(t))\le x-c 
\big\} , $$
$$ \P\big\{\inf_{a\le t\le b}|W(t)-M|\ge  c\,\big|\,W(a) =M-x\big\}=\P\big\{\sup_{a\le t\le b} (W(t)-W(a))\le x-c  \big\} .
$$
As   for $y\ge0$, (\cite{CR}, Theorem 1.5.1) 
$$ \P\{ \sup_{0\le t\le T} W(t)>y\}= 2\P\{W(T)>y\},$$
we get  if $|x|>c$, 
\begin{equation} \P\big\{\inf_{a\le t\le b}|W(t)-M|\ge  c\,\big|\, W(a)=M+ x\big\}     =1-2\Psi({|x|-c\over \sqrt{b-a}})  .
\end{equation} 
Therefore
\begin{eqnarray}&{}&\P\big\{\inf_{a\le t\le b}|W(t)-M|\ge  c \big\}\cr&=&\int_\R\P\big\{\inf_{a\le t\le b}|W(t)-M|\ge  c\,\big|\,W(a)
=u\big\}{ e^{-{u^2\over 2a}}\over
\sqrt{2\pi a}}du
\cr&=&\int_{|u-M|>c}\P\big\{\inf_{a\le t\le b}|W(t)-M|\ge  c\,\big|\,W(a) =u\big\}{
e^{-{u^2\over 2a}}\over
\sqrt{2\pi a}}du
 \cr &  =& \int_{|u-M|>c } \Big[1-2\Psi({|u-M|-c\over \sqrt{b-a}})\Big]{ e^{-{u^2\over
2a}}\over
\sqrt{2\pi a}}du\cr &=&\int_{|v|>c } \Big[1-2\Psi({|v|-c\over \sqrt{b-a}})\Big]{ e^{-{(M+v)^2\over
2a}}\over
\sqrt{2\pi a}}dv   ,\end{eqnarray}
as claimed.
  \cqfd
 
\begin{rem} \rm It follows from Lemma \ref{l2} 
  that 
\begin{eqnarray} \P\big\{\inf_{a\le t\le b}|W(t)-M|>  0 \big\}&=&\lim_{c\downarrow 0}\P\big\{\inf_{a\le t\le b}|W(t)-M|\ge  c \big\}\cr &
=&
\int_{\R} \Big[1-2\Psi({|v| \over \sqrt{b-a}})\Big]{ e^{-{(M+v)^2\over
2a}}\over
\sqrt{2\pi a}}dv 
   .\end{eqnarray}
Thus 
 $$\P\big\{\inf_{a\le t\le b}|W(t)-M|\ge 0 \big\}=1\not=\int_{\R} \Big[1-2\Psi({|v| \over \sqrt{b-a}})\Big]{ e^{-{(M+v)^2\over
2a}}\over
\sqrt{2\pi a}}dv 
, $$
     yielding a discontinuity at $0$. 
 Take for instance $a=1$, $b=1+\m^2$; the integral above is 
 $$\int_{\R} \Big[1-2\Psi({|v| \over \m})\Big]{ e^{-{(M+v)^2\over
2 }}\over
\sqrt{2\pi }}dv 
\ \to \  0 \qq  \qq \m\to\infty.$$

And 
\begin{equation}\P\big\{\inf_{a\le t\le b}|W(t)-M|= 0
\big\}
 =2\int_{\R}    \Psi({|v| \over \sqrt{b-a}}) { e^{-{(M+v)^2\over
2a}}\over
\sqrt{2\pi a}}dv.
\end{equation}
  \end{rem}
\smallskip We will also show

\begin{lem} \label{l3} There exists an  absolute constant $C$, such that for every real $M$, and $0<a<b<\infty$, 
\begin{eqnarray*}\P\big\{\inf_{a\le t\le b}|W(t)-M|= 0 \big\} &=&2\int_{\R}    \Psi({|v| \over \sqrt{b-a}}) { e^{-{(M+v)^2\over
2a}}\over
\sqrt{2\pi a}}dv \cr &\le&   C \min\Big(1,\sqrt{b-a \over
     a}  e^{-{ M ^2\over
8\max(a, b-a)}} \Big)   .
\end{eqnarray*}
In particular $$\P\big\{\inf_{a\le t\le b}|W(t)-M|= 0 \big\}\le C\min\Big(1,\sqrt{{b-a\over  a}}\Big).$$
\end{lem}{\it Proof.} By (2.6) 
\begin{eqnarray}&{}&\P\big\{\inf_{a\le t\le b}|W(t)-M|= 0
\big\} \cr&=&2\int_{\R}    \Psi({|v| \over \sqrt{b-a}}) { e^{-{(M+v)^2\over
2a}}\over
\sqrt{2\pi a}}dv\cr &=&2\bigg\{\int_{|v|\le {M\over 2}}  + \int_{|v|> {M\over 2}} \bigg\}   \Psi({|v| \over \sqrt{b-a}}) {
e^{-{(M+v)^2\over 2a}}\over
\sqrt{2\pi a}}dv \cr &\le&\sqrt{{2 \over
  \pi a}}\Bigg\{ e^{-{ M ^2\over
8a}}\int_{\R}    \Psi({|v| \over \sqrt{b-a}})  dv  + \int_{|v|> {M\over 2}}    \Psi({|v| \over \sqrt{b-a}})  dv\Bigg\}
\cr &=&\sqrt{{2 \over
  \pi a}}\Bigg\{ e^{-{ M ^2\over
8a}}\sqrt{b-a}\int_{\R}    \Psi( |w| )  dw   + \int_{|v|> {M\over 2}}    \Psi({|v| \over \sqrt{b-a}})  dv\Bigg\}  . 
\end{eqnarray}
    Recall that
  the Mills'  ratio $\displaystyle{R(x) =
 e^{x^2/2}\int_x^\infty e^{-t^2/2}\ dt=(\sqrt{2\pi}})e^{x^2/2}\Psi(x) $  verifies for all $x\ge 0$ 
 $  R(x) \le 
 \sqrt{{\pi/ 2}}$. 
 Thus
 \begin{eqnarray} \int_{|v|> {M\over 2}}    \Psi({|v| \over \sqrt{b-a}})  dv &=&\sqrt{b-a} \int_{|w|> {M\over 2\sqrt{b-a}}}    \Psi(|w|) 
dw\cr & \le & C\sqrt{b-a}
\int_{|w|> {M\over 2\sqrt{b-a}}}     e^{-w^2/2} dw \cr & \le & C\sqrt{b-a}      e^{-{M^2\over 8(b-a)}}.\end{eqnarray}
 Therefore
\begin{eqnarray}\P\big\{\inf_{a\le t\le b}|W(t)-M|= 0 \big\}&\le&  C \sqrt{b-a \over
     a} \Big\{ e^{-{ M ^2\over
8a}}   + e^{-{M^2\over 8(b-a)}} \Big\}\cr &\le & C \sqrt{b-a \over
     a}  e^{-{ M ^2\over
8\max(a, b-a)}}    . \end{eqnarray}
 \cqfd 
 
One can recover as a special case that (see \cite{D} p.248)
$$\P\big\{\inf_{a\le t\le b}|W(t)|= 0 \big\}  = 1- { 2\over   \pi  }\arctan  \sqrt{ {   a  \over b-a}} ,$$
or, equivalently
 $\P\big\{  W(t) \, \hbox{has no zero in $(a,b)$} \big\}  =   ({ 2/   \pi  })\arcsin \sqrt{ {   a / b }}  $. 
 It is possible to also give an exact expression of the
probability $\P\big\{\inf_{a\le t\le b}|W(t)-M|= 0 \big\}$, although for $M\not= 0$ this one is relatively more complicated. This is indicated in the Lemma
below.

\begin{lem} \label{l4}  Let $0<a\le b<\infty$. We have 
$$\P\big\{\inf_{a\le t\le b}|W(t)-M|= 0 \big\}    =  -\sqrt{2\over
\pi}   {M\over \sqrt a}  \int_{   \sqrt{ {b-a\over b}  }}^1  u e^{-{ (M u)^2\over 2a}      } 
    \Big[\int_{ |x|\le M  \sqrt{{ 1-u^2}\over a} }    {e^{-{x^2 \over 2}  }dx\over \sqrt{2\pi}} \Big]    du 
$$$$ \quad +\Big[1-\big({2\over
\pi}\big) e^{-{  M  ^2\over 2a}   }\arctan \sqrt{{a\over b-a}}\Big] . 
$$ 

 In particular    
 $\P\big\{\inf_{a\le t\le b}|W(t)|= 0 \big\}  = 1- { 2\over   \pi  }\arctan  \sqrt{ {   a  /( b-a)}}$.   
 And
for every positive real
$c$  
\begin{eqnarray*}\P\big\{0<\inf_{a\le t\le b}|W(t)|<c \big\}& = & 2 \int_{0}^{c/\sqrt a} \big(1-2\Psi(  
  { u\sqrt a  \over \sqrt{b-a}}
 )\big){ e^{-{u^2\over 2 }}\over \sqrt{2\pi }}du\cr &{}&\ +  4\int_{c/\sqrt a}^\infty \Big( \Psi(  
  { u\sqrt a-c \over \sqrt{b-a}}
 )-  \Psi(  
  { u\sqrt a  \over \sqrt{b-a}}
 )\Big){ e^{-{u^2\over 2 }}\over \sqrt{2\pi }}du  .
\end{eqnarray*} 
\end{lem}

\noi{\it Proof.}  Let $M_1=M/\sqrt{  a}$, and put   $F(s)= \int_{\R}   2\Psi\big( |w|s\big)   e^{-{(M_1+w )^2\over
2 }}{dw\over
\sqrt{2\pi  }} $.  We have
\begin{eqnarray}\P\big\{\inf_{a\le t\le b}|W(t)-M|= 0 \big\}& =&\int_{\R}   2\Psi({|v| \over \sqrt{b-a}}) { e^{-{(M+v)^2\over
2a}}\over
\sqrt{2\pi a}}dv\cr &=&\int_{\R}   2\Psi\Big( |w|\sqrt{a \over  {b-a}}\Big) { e^{-{(M_1+w )^2\over
2 }}\over
\sqrt{2\pi  }}dw \cr &=&F\Big(\sqrt{a \over  {b-a}}\Big). 
\end{eqnarray}

As ${\partial\over \partial
s}\Psi(|w|s)=|w|\Psi'(|w|s)=-{|w|\over \sqrt{2\pi}}e^{-(|w|s)^2/2}$, we have
\begin{eqnarray*}{\partial\over \partial s}F(s)&= &\int_{\R}   2{\partial\over \partial s}\Psi\Big((|w|s\Big) { e^{-{(M_1+w )^2\over
2 }}\over
\sqrt{2\pi  }}dw \cr & =&-\sqrt{{2\over
 \pi  }} \int_{\R}    |w| e^{-{1\over 2}[ ( w s)^2  + (M_1+w )^2 
 ] }{dw \over
\sqrt{2\pi  }}. 
\end{eqnarray*}
But $[ ( w s)^2  + (M_1+w )^2 
 ]=[w\sqrt{s^2+1}+ {M_1 \over \sqrt{s^2+1}}]^2+ {M_1^2s ^2\over s^2+1} $,
 hence
\begin{eqnarray}{\partial\over \partial s}F(s)& =&-\sqrt{{2\over
 \pi  }}e^{-   { M_1^2s ^2\over 2(s^2+1)} 
   } \int_{\R}    |w| e^{-{1\over 2} [w\sqrt{s^2+1}+ {M_1 \over \sqrt{s^2+1}}]^2   }{dw \over
\sqrt{2\pi  }} \cr &
=&-\sqrt{2\over
\pi} \ {e^{-   { M_1^2s ^2\over 2(s^2+1)} 
   } \over s^2+1}\int_{\R}    |z| e^{-{1\over 2} [z+ {M_1 \over \sqrt{s^2+1}}]^2   }{dz \over
\sqrt{2\pi  }  }\cr &=&-\sqrt{2\over
\pi} \ {e^{-   { M_1^2s ^2\over 2(s^2+1)} 
   } \over s^2+1}\ \E \Big|g-{M_1 \over \sqrt{s^2+1}}\Big|, \end{eqnarray} 
where $g$ denotes a Gaussian standard random variable. But  for any real 
$a$
\begin{equation}\E |g+a|=|a|\int_{-|a|}^{|a|} e^{-x^2/2} {dx\over \sqrt{2\pi}}+ \sqrt{2\over
\pi} \ e^{-a^2/2}
 .
\end{equation}
Hence
\begin{eqnarray}{\partial\over \partial s}F(s)&=&-\sqrt{2\over
\pi} \, {e^{-   { M_1^2s ^2\over 2(s^2+1)} 
   } \over s^2+1}\bigg\{{M_1 \over \sqrt{s^2+1}}\int_{ |x|\le {M_1 \over \sqrt{s^2+1}}}    {e^{-x^2/2}dx\over \sqrt{2\pi}}+
\sqrt{2\over
\pi} \, e^{-{ M_1 ^2\over 2( s^2+1) }  }\bigg\}\cr &  =&-\sqrt{2\over
\pi}     {M_1   e^{-   { M_1^2s ^2\over 2(s^2+1)} 
   } \over(s^2+1)^{3/2}}\int_{ |x|\le {M_1 \over \sqrt{s^2+1}}}   {e^{-x^2/2}dx\over \sqrt{2\pi}} 
- {2\over
\pi}\, {e^{-  M_1^2/2   } \over s^2+1}  .  \end{eqnarray}
 
If $M=0$, this takes a much simplified form
$${\partial\over \partial s}F(s) = - {2\over
\pi} \, {1 \over s^2+1}. $$
  Further $F(0)=  \int_{\R}       e^{-{(M_1+w )^2\over
2 }}{dw\over
\sqrt{2\pi  }}=1 $. Therefore
\begin{eqnarray}  F(s)- 1   & = &-\sqrt{2\over
\pi} M_1\int_0^s  e^{-   { M_1^2t ^2\over 2(t^2+1)}}  \Big[\int_{ |x|\le {M_1 \over \sqrt{t^2+1}}}  e^{- x^2/2 } {dx\over \sqrt{2\pi}} \Big] { 
 dt\over (t^2+1)^{3/2}} \cr &{}& - {2\over
\pi}\, e^{- M_1^2 /2  }\int_0^s{dt \over t^2+1} 
\cr &=& - e^{- M_1^2 /2  }\bigg\{\sqrt{2\over
\pi} M_1\int_0^s    \Big[\int_{ |x|\le {M_1 \over \sqrt{t^2+1}}}  e^{- x^2/2 } {dx\over \sqrt{2\pi}} \Big] { 
 e^{    { M_1^2 \over 2(t^2+1)}} dt\over (t^2+1)^{3/2}} \cr &{}& - {2\over
\pi}  \arctan s   \bigg\}.
\end{eqnarray} 
  
Consequently, 
\begin{eqnarray}&{}&\P\big\{\inf_{a\le t\le b}|W(t)-M|= 0 \big\}   =F\Big(\sqrt{a \over  {b-a}}\Big)  
\cr & =& -\sqrt{2\over
\pi} {M\over \sqrt a}\int_{   \sqrt{ {b-a\over b}  }}^1  u e^{-{ (M u)^2\over 2a}      } 
    \Big[\int_{ |x|\le M  \sqrt{{ 1-u^2}\over a} }  e^{-{x^2 \over 2}  } {dx\over \sqrt{2\pi}} \Big]    du 
 \cr &{}&  +\Big[1- {2\over
\pi} \, e^{-{  M  ^2\over 2a}   }\arctan \sqrt{{a\over b-a}}\Big] . \end{eqnarray}
  If $M =0$, this is simplified into
\begin{eqnarray}\P\big\{\inf_{a\le t\le b}|W(t) |= 0 \big\}  & =\Big[1- {2\over
\pi}\,  \arctan \sqrt{{a\over b-a}}\Big] . \end{eqnarray}

  As concerning the second formula,  it suffices to write
\begin{eqnarray*} \P\big\{0<\inf_{a\le t\le b}|W(t)|\le c
\big\}  &=&\P\big\{
\inf_{a\le t\le b}|W(t)|\le c \big\}-\P\big\{ \inf_{a\le t\le b}|W(t)|=0
\big\}\cr &=&\Big(1-2\int_{c\over \sqrt a}^\infty \Big(1-2\Psi(  
  { u\sqrt a-c \over \sqrt{b-a}}
 )\Big){ e^{-{u^2\over 2 }}\over \sqrt{2\pi }}du\Big)\cr&{}&-\Big(1-2\int_{0}^\infty \big(1-2\Psi(  
  { u\sqrt a  \over \sqrt{b-a}}
 )\big){ e^{-{u^2\over 2 }}\over \sqrt{2\pi }}du\Big)\cr&=&-   2\int_{c\over\sqrt a}^\infty \Big(1-2\Psi(  
  { u\sqrt a-c \over \sqrt{b-a}}
 )\Big){ e^{-{u^2\over 2 }}\over \sqrt{2\pi }}du \cr&{}&+  2\int_{0}^\infty \big(1-2\Psi(  
  { u\sqrt a  \over \sqrt{b-a}}
 )\big){ e^{-{u^2\over 2 }}\over \sqrt{2\pi }}du \cr&=&  2\int_{0}^{c\over \sqrt a} \big(1-2\Psi(  
  { u\sqrt a  \over \sqrt{b-a}}
 )\big){ e^{-{u^2\over 2 }}\over \sqrt{2\pi }}du \cr&{ }&\!\!+    4\int_{c\over \sqrt a}^\infty \Big[ \Psi(  
  { u\sqrt a-c \over \sqrt{b-a}}
 )-  \Psi(  
  { u\sqrt a  \over \sqrt{b-a}}
 )\Big){ e^{-{u^2\over 2 }}du\over \sqrt{2\pi }} \Big]    .
\end{eqnarray*}
\cqfd

Notice that $    {\pi\over 2}=\int_0^\infty { dt\over     (1+t^2)}$, and so
\begin{eqnarray*}1- {2\over
\pi}\,  \arctan s&=& {2\over
\pi}\,\big[\int_0^\infty { dt\over     (1+t^2)}-\int_0^s { dt\over     (1+t^2)}\big]\cr &= & {2\over
\pi}\, \int_s^\infty { dt\over     (1+t^2)} \le C\min(1,s^{-1}).\end{eqnarray*}
  We deduce the bound obtained in Lemma 2.2
 $$\P\big\{\inf_{a\le t\le b}|W(t)|= 0 \big\}  \le  C\min\Big(1,\sqrt{{b-a\over  a}}\Big). $$ 
 
 \bigskip

In what follows we shall be interested in finding estimates of the delicate random variable $  \b_{[a,b]}^{ -\a} \cdot\chi_{\{\b_{[a,b]}> 0\}}$,   $0\le
\a<1$, where we set
$$\b_{[a,b]}  :=\inf_{a\le t\le b}|W(t)|. $$

\begin{prop} Let $b>a>0$, and $\eta>0$. Then,  
$$ \P  \big\{0<\b_{[a,b]}\le \eta\big\} \le     {  16\eta \over  
\sqrt{ 2\pi   }} \min\Big(   { 1  
\over  \sqrt{ b-a }}      \, ,{1\over \sqrt{    a} }\Big) .  $$
Further, for any real $\a$, $0\le\a<1$   $$   \E \ \Big\{{1\over \b_{[a,b]}^{ \a}}\cdot\chi\{\b_{[a,b]}> 0\}
\Big\}\le  
 {47\over  1-\a   }\min\Big(   { 1  
\over  \sqrt{ b-a }}      \, ,{1\over \sqrt{  a} }\Big)  +1.   $$
\end{prop}
\noi The result gives a   control  which is uniform in $b-a$.   We have $\b_{[a,b]}\to 0$ as 
$b-a\to \infty$, but in the same time   the constraint $\b_{[a,b]}> 0$ becomes also stronger, making the
probability of the set
$\{\b_{[a,b]}> 0\}$   small. When $b-a\to 0$,  $\b_{[a,b]}\to |W(a)|$, and this is reflected by the term ${1/ \sqrt   a  }$
in our estimate.
\medskip\par\noi{\it Proof.} Write $\b =\b_{[a,b]}   $.  By using the second formula in Lemma \ref{l4}   
 \begin{eqnarray}\P  \big\{0<\b\le \eta\big\}& =&  4 \int_{0}^{\eta/  \sqrt a } \big(1-2\Psi(  
  { u\sqrt a  \over \sqrt{b-a}}
 )\big){ e^{-{u^2\over 2 }}\over \sqrt{2\pi }}du\cr &{}& +  4\int_{\eta/  \sqrt a }^\infty \Big( \Psi(  
  { u\sqrt a-\eta\over \sqrt{b-a}}
 )-  \Psi(  
  { u\sqrt a  \over \sqrt{b-a}}
 )\Big){ e^{-{u^2\over 2 }}\over \sqrt{2\pi }}du . \end{eqnarray}
As $1-2\Psi(x)=\int_{-x}^x e^{-u^2/2} {du \over \sqrt
{2\pi}} \le \min(  (2/\pi)^{1/2}x,1)$, $x\ge 0$
\begin{eqnarray}  &{}&  \int_{0}^{\eta/  \sqrt a} \big(1-2\Psi(  
  { u\sqrt a  \over \sqrt{b-a}}
 )\big){ e^{-{u^2\over 2 }}\over \sqrt{2\pi }}du \cr
&\le  &      \int_{0}^{\eta/  \sqrt a}\min\Big(  \big({2\over \pi}\big)^{1/2}   { u
\sqrt a 
\over \sqrt{b-a}},1\Big)  { e^{-{u^2\over
2 }}\over \sqrt{2\pi }}du   
  \cr
&\le  &  \min\Big(  \big({2\over \pi}\big)^{1/2}   { \sqrt a 
\over \sqrt{b-a}} \int_{0}^{\eta/  \sqrt a}   u { e^{-{u^2\over
2 }}\over \sqrt{2\pi }}du\, ,{\eta\over  \sqrt{ 2\pi a} }\Big)   
  \cr &\le &   \min\Big(   { \eta \max_{u\ge 0} u   e^{-{u^2\over
2 }}  
\over  \pi\sqrt{ b-a }}      \, ,{\eta\over  \sqrt{ 2\pi a} }\Big)   
  =\eta\min\Big(   { 1  
\over \pi\sqrt{ b-a }}      \, ,{1\over \sqrt{ 2\pi a} }\Big) .
\end{eqnarray}
Furthermore
\begin{equation} \Psi(  
  { u\sqrt a-\eta \over \sqrt{b-a}}
 )-  \Psi(  
  { u\sqrt a  \over \sqrt{b-a}}
 )=\int_{{ u\sqrt a-\eta \over \sqrt{b-a}}
}^{{ u\sqrt a  \over \sqrt{b-a}}}  e^{-{u^2\over 2 }} { du\over \sqrt{2\pi }} \le    { \eta \over  \sqrt{2\pi (b-a)}}, 
\end{equation}
which implies
\begin{eqnarray}\int_{\eta/ \sqrt a }^\infty \Big( \Psi(  
  { u\sqrt a-\eta \over \sqrt{b-a}}
 )-  \Psi(  
  { u\sqrt a  \over \sqrt{b-a}}
 )\Big){ e^{-{u^2\over 2 }}\over \sqrt{2\pi }}du&\le&     { \eta \over  \sqrt{2\pi (b-a)}} \int_{\eta/ \sqrt a }^\infty  {
e^{-{u^2\over 2 }}\over \sqrt{2\pi }}du\cr & \le  &   { \eta \over 
\sqrt{2\pi (b-a)}}  .\end{eqnarray}
Besides, with the variable  change $u= v\sqrt{(b-a)\over a}$, letting $v_\eta=    { \eta \over  
\sqrt{   b-a }}$.
$$\int_{{\eta/ \sqrt a }}^\infty \Big( \Psi(  
  { u\sqrt a-\eta \over \sqrt{b-a}}
 )-  \Psi(  
  { u\sqrt a  \over \sqrt{b-a}}
 )\Big){ e^{-{u^2\over 2 }}\over \sqrt{2\pi }}du \qq\qq\qq\qq$$$$\qq\qq\qq =\sqrt{ b-a \over a} \int_{v_\eta}^\infty \big( \Psi(  
  v-v_\eta
 )-  \Psi(v)\big){ e^{-{v^2(b-a)\over 2a }}\over \sqrt{2\pi }}dv. $$
We have $$\sqrt{ b-a \over a} \int_{v_\eta}^{2v_\eta} \big( \Psi(  
  v-v_\eta
 )-  \Psi(v)\big){ e^{-{v^2(b-a)\over 2a }}\over \sqrt{2\pi }}dv\le  \sqrt{ b-a \over a}{v_\eta\over \sqrt{2\pi }}=   { \eta \over  \sqrt{
2\pi  a }} .
 $$
Now if $v\ge 2v_\eta$, then $v-v_\eta\ge v-v/2=v/2$. Consequently
\begin{eqnarray*}&{}&\sqrt{ b-a \over a} \int_{2v_\eta}^\infty \big( \Psi(  
  v-v_\eta
 )-  \Psi(v)\big){ e^{-{v^2(b-a)\over 2a }}\over \sqrt{2\pi }}dv\cr &=&\sqrt{ b-a \over a} \int_{2v_\eta}^\infty\Big\{
\int_{v-v_\eta}^v{e^{-{x^2
\over 2  }}\over \sqrt{2\pi }}dx\Big\} { e^{-{v^2(b-a)\over 2a }}\over \sqrt{2\pi }}dv\cr &\le &\sqrt{ b-a \over a}v_\eta
\int_{2v_\eta}^\infty { e^{-{ v ^2 / 8  }} \over \sqrt{2\pi }}  { e^{-{v^2(b-a)\over 2a }}\over \sqrt{2\pi }}dv\cr &\le &\sqrt{ b-a \over
a}v_\eta \int_{0}^\infty  e^{-{ v ^2 / 8  }}    { dv\over 2\pi  } ={ \eta \over  
\sqrt{ 2\pi  a }} .
\end{eqnarray*}
It follows that
\begin{eqnarray}\int_{{\eta/ \sqrt a }}^\infty \Big( \Psi(  
  { u\sqrt a-\eta \over \sqrt{b-a}}
 )-  \Psi(  
  { u\sqrt a  \over \sqrt{b-a}}
 )\Big){ e^{-{u^2\over 2 }}\over \sqrt{2\pi }}du&  \le \min\Big({ 2\eta \over  
\sqrt{ 2\pi  a }},   { \eta \over   
\sqrt{2\pi (b-a)}}  \Big) 
\cr &\le  { 2\eta \over  
\sqrt{ 2\pi   }} \min\Big({ 1 \over 
\sqrt{    a }},   { 1 \over 
\sqrt{  b-a }}  \Big) .\end{eqnarray}

By reporting 
 \begin{eqnarray}\P  \big\{0<\b\le \eta\big\}& \le      4\eta \min\Big(   { 1  
\over \pi\sqrt{ b-a }}      \, ,{1\over \sqrt{ 2\pi a} }\Big)  +    {  8\eta \over  
\sqrt{ 2\pi   }} \min\Big({ 1 \over 
\sqrt{    a }},   { 1 \over 
\sqrt{  b-a }}  \Big) \cr &\le     {  16\eta \over  
\sqrt{ 2\pi   }} \min\Big(   { 1  
\over  \sqrt{ b-a }}      \, ,{1\over \sqrt{    a} }\Big) . \end{eqnarray}

Now  let
$X$ be a random variable such that
$X\ge 0$ a.s. As
$\chi\{ ]0,1]\}=\sum_{n=0}^\infty
\chi\{]{1\over 2^{n+1}} , {1\over 2^{n }}]\}$, we have by the Beppo-Levi Theorem
\begin{equation}\E {1\over X^{ \a}}\cdot\chi\{X> 0\}= \sum_{n=0}^\infty \E \Big({1\over X^{ \a}}\cdot\chi\{{1\over 2^{n+1}}<X\le {1\over
2^{n }}\}\Big)+ \E {1\over X^{ \a}}\cdot\chi\{X> 1\}, 
\end{equation}
where "$=$" means   that both expresions are simultaneously finite or infinite. And
\begin{eqnarray} &{}  &\sum_{n=0}^\infty \E \Big({1\over X^{ \a}}\cdot\chi\{{1\over 2^{n+1}}<X\le {1\over
2^{n }}\}\Big)+ \E {1\over X^{ \a}}\cdot\chi\{X> 1\}\cr &\le &\sum_{n=0}^\infty 2^{\a(n+1)}\P  \big\{0<X\le {1\over 2^{n }}\big\}  +
1.
 \end{eqnarray}

Apply this with $X=\b$. Then $\E  \b^{ -\a} \cdot\chi\{\b> 0\}$ will be finite once we prove that the series
$$\sum_{n=0}^\infty 2^{\a(n+1)}\P  \big\{0<\b\le {1\over 2^{n }}\big\}$$
is convergent. By using (2.13) with $\eta=2^{-n}$, we get
\begin{eqnarray}\sum_{n=0}^\infty 2^{\a(n+1)}\P  \big\{0<\b\le {1\over 2^{n }}\big\}&\le&  {   2^{\a +4} \over  
\sqrt{ 2\pi   }}  \min\Big(   { 1  
\over  \sqrt{ b-a }}      \, ,{1\over \sqrt{  a} }\Big)\sum_{n=0}^\infty   2^{-n(1-\a)} \cr 
&\le & {32\over 2^{   1-\a }-1  }\min\Big(   { 1  
\over  \sqrt{ b-a }}      \, ,{1\over \sqrt{  a} }\Big)\cr &\le & {47\over  (  1-\a)   }\min\Big(   { 1  
\over  \sqrt{ b-a }}      \, ,{1\over \sqrt{  a} }\Big)    .\end{eqnarray}
   And we conclude that
 $$   \E {1\over \b^{ \a}}\cdot\chi\{\b> 0\} \le  
 {47\over  1-\a   }\min\Big(   { 1  
\over  \sqrt{ b-a }}      \, ,{1\over \sqrt{  a} }\Big)+1 ,   $$
as claimed. \cqfd
 
\section{Proof of Theorem} Recall that $I_N=[T_N,T_{N+1}]$. It is now easy. As $$\sum_{N\ge 1} \P 
\big\{0<\b_{I_N}\le \eta_N\big\}
\le  C\sum_{N\ge 1}\eta_N \min\Big(   { 1  
\over  \sqrt{ T_{N+1}-T_{N } }}      \, ,{1\over \sqrt{   T_{N }} }\Big)<\infty ,  $$
we deduce from Borel-Cantelli Lemma that
$$\P\Big\{ \inf_{t\in I_N}|W(t)|\ge  \eta_N\ {\rm  or}\   W(t) =0 \ {\rm  for \ some}\ t\in I_N, \quad  N\ {\rm
ultimately}\Big\}=1.
$$
\cqfd 
  
 \noi  {\sl Acknowledgments.} {\it I thank Istvan Berkes for valuable remarks and careful reading of the paper.}
 
%%%%%%%%%%%%%%%%%%%%%%%%%%%%%%%%%%%%%%%%%%%%%%%%%%%%%%%%%%%%%%%%%%%%%%%%%%%%%%%%%%%%%%%%%%%%%%%%
%%%%%%%%%%%%%%%%%%%%%%%%%%%%%%%%%%%%%%%%%%%%%%%%%%%%%%%%%%%%%%%%%%%%%%%%%%%%%%%%%%%%%%%%%%%%%%%%
%%%%%%%%%%%%%%%%%%%%%%%%%%%%%%%%%%%%%%%%%%%%%%%%%%%%%%%%%%%%%%%%%%%%%%%%%%%%%%%%%%%%%%%%%%%%%%%%
%%%%%%%%%%%%%%%%%%%%%%%%%%%%%%%%%%%%%%%%%%%%%%%%%%%%%%%%%%%%%%%%%%%%%%%%%%%%%%%%%%%%%%%%%%%%%%%%
%%%%%%%%%%%%%%%%%%%%%%%%%%%%%%%%%%%%%%%%%%%%%%%%%%%%%%%%%%%%%%%%%%%%%%%%%%%%%%%%%%%%%%%%%%%%%%%%
%%%%%%%%%%%%%%%%%%%%%%%%%%%%%%%%%%%%%%%%%%%%%%%%%%%%%%%%%%%%%%%%%%%%%%%%%%%%%%%%%%%%%%%%%%%%%%%%
%%%%%%%%%%%%%%%%%%%%%%%%%%%%%%%%%%%%%%%%%%%%%%%%%%%%%%%%%%%%%%%%%%%%%%%%%%%%%%%%%%%%%%%%%%%%%%%%
%%%%%%%%%%%%%%%%%%%%%%%%%%%%%%%%%%%%%%%%%%%%%%%%%%%%%%%%%%%%%%%%%%%%%%%%%%%%%%%%%%%%%%%%%%%%%%%%

{
}

 \noi {\phh Michel  Weber, \noi  Math\'ematique (IRMA),
Universit\'e Louis-Pasteur et C.N.R.S.,   7  rue Ren\'e Descartes,
67084 Strasbourg Cedex, France.
\par\noindent
E-mail: \  \tt weber@math.u-strasbg.fr} 

\baselineskip 12pt
\begin{thebibliography}{99}
 \bibitem{CR}  \vskip 1pt\noi   Cs\"org\"o M., R\'ev\'esz P.  [1981] {\sl Strong approximations in probability and Statistics},
Akad\'emiai Kiad\'o, Budapest.
 \bibitem{D}  \vskip 1pt\noi   Durrett R.  [1996] {\sl Stochastic calculus: A practical introduction},
Prob.\& Stoch.
\end{thebibliography}
\end{document}